\documentclass[11pt,reqno]{amsart}
\usepackage{graphicx}
\usepackage{amssymb,latexsym,epsfig,mathrsfs,graphpap,graphics,amsmath,amssymb}
\usepackage{amstext}
\usepackage{verbatim}
\usepackage{color}

\newtheorem{conj.}[thm]{Conjecture}

\theoremstyle{definition}

\theoremstyle{remark}

\numberwithin{equation}{section}

\begin{document}

\begin{flushleft}
  {\bf\Large {Uncertainty Principles for the Continuous Shearlet  \\[1.6mm]Transforms in Arbitrary Space Dimensions}}
\end{flushleft}

\parindent=0mm \vspace{.4in}

  {\bf{Firdous A. Shah$^{\star}$  and   Azhar Y. Tantary$^{\star}$ }}

 \parindent=0mm \vspace{.1in}
{\small \it $^{\star}$Department of  Mathematics,  University of Kashmir, South Campus, Anantnag-192101, Jammu and Kashmir, India. E-mail: $\text{fashah@uok.edu.in}$;\,$\text{aytku92@gmail.com}$}

\parindent=0mm \vspace{.2in}
{\small {\bf Abstract.} The aim of this article is to formulate some novel uncertainty principles for the continuous shearlet transforms in arbitrary space dimensions. Firstly, we derive an analogue of the Pitt's inequality for the continuous shearlet transforms, then we formulate the Beckner's uncertainty principle via two approaches: one based on a sharp estimate from Pitt's inequality and the other from the classical Beckner's inequality in the Fourier domain.  Secondly, we consider a logarithmic Sobolev inequality for the continuous shearlet transforms which has a dual relation with Beckner's inequality. Thirdly, we derive Nazarov's uncertainty principle for the shearlet transforms which shows that it is impossible for a non-trivial function and its shearlet transform to be both supported on sets of finite measure. Towards the culmination, we formulate local uncertainty principles for the continuous shearlet transforms in arbitrary space dimensions.

\parindent=0mm \vspace{.1in}
{\bf{Keywords:}} Shearlets. Uncertainty principle.  Pitt’s inequality. Beckner's inequality. Sobolev inequality. Nazarov's uncertainty principle. Local uncertainty principle. Fourier transform.

\parindent=0mm \vspace{.1in}
{\bf {Mathematics Subject Classification:}} 26D10. 35A23. 42B10. 42C40. 42A38.}

\section{Introduction}

\parindent=0mm \vspace{.0in}
Shearlets are the outcome of a series of  multiscale methods such as wavelets, ridgelets, curvelets, contourlets and many others introduced during the last few decades  with the aim to achieve optimally sparse approximations for higher dimensional signals by employing the basis elements with much higher directional sensitivity and various shapes \cite{Lab,Kuty,Dal1,Dal2}. Unlike the classical wavelets, shearlets are non-isotropic in nature, they offer optimally sparse representations, they
allow compactly supported analyzing elements, they are associated with fast decomposition algorithms and they provide a unified treatment of continuum and digital data. However, similar to the wavelets, they are an affine-like system of well-localized waveforms at various scales, locations and orientations; that is, they are generated by dilating and translating one single generating function, where the dilation matrix is the product of a parabolic scaling matrix and a shear matrix and hence, they are a specific type of composite dilation wavelets \cite{JJ,SZS,ABHK,DS1,DS2}. The importance of shearlet transforms have been  widely acknowledged and since their inception, they have emerged as one of the most effective frameworks for representing multidimensional data ranging over the areas of signal and image processing, remote sensing, data compression, and several others, where the detection of directional structure of the analyzed signals play a role \cite{GL1,GL2}.

\parindent=8mm \vspace{.1in}
For any  $f\in L^2( {\mathbb R}^{n})$,  the continuous shearlet transform in arbitrary space dimension is defined by \cite{Dal2}
\begin{align*}
{\mathcal {SH}}_{\psi}f(a,s,t)=\Big\langle{f,\psi_{a,s,t}}\Big\rangle=\int_{\mathbb R^{n}} f(x)\,\overline{\psi_{a,s,t}(x)}\,dx, \tag{1.1}
\end{align*}
where $\psi_{a,s,t}(x)=|\det A_{a}|^{\frac{1}{2n}-1}\psi\big(A_{a}^{-1}S_{s}^{-1} \left(x-t\right)\big),a\in\mathbb R\setminus \left\{0\right\}, s\in\mathbb R^{n-1}, t\in\mathbb R^n$ is the shearlet family constituted by the combined action of the scaling $D_{A_{a}}$, sharing $\mathcal{D}_{S_{s}}$ and translation ${T_{t}}$  operators on the analyzing function $\psi \in L^2(\mathbb R^n)$ given by
 \begin{align*}
D_{A_{a}}\psi(x)=|\det A_{a}|^{-1/2}\,{\psi}\left({A_{a}^{-1}}x\right),~~\mathcal{D}_{S_{s}}\psi (x)={\psi}\left({S_{s}^{-1}}x\right),~~\text{and}~~  {T_{t}}\psi(x)={\psi}(x-t), \tag{1.2}
\end{align*}
 respectively, and the matrices involved in (1.2) are given by
\begin{align*}
A_{a}=\left(\begin{array}{cc}
 a & {\bf 0}_{n-1}^{T} \\{\bf 0}_{n-1} & \text{sgn}(a)\, a^{1/n} \,I_{n-1}\\
 \end{array}\right) ~~\text{ and}~~
 S_{s}=\left(\begin{array}{cc} 1 & {\bf s}^{T} \\ {\bf 0}_{n-1} & I_{n-1}
\end{array} \right),\tag{1.3}
\end{align*}
 ${\bf s}^{T}=\big(s_1, s_2,\dots, s_{n-1}\big)$, $\text{sgn}(\cdot)$ and ${\bf 0}$  denotes the well known  Signum function and the  null vector, respectively. For the brevity, we shall rewrite the shearlet family $\psi_{a,s,t}(x)$ as
  \begin{align*}
\psi_{a,s,t}(x)=\big|\det M_{sa}\big|^{-1/2}\psi\Big({ M_{sa}}^{-1} \left(x-t\right)\Big),\tag{1.4}
\end{align*}
where $M_{sa}=S_{s}A_{a}$ is the composition of the  parabolic scaling matrix $A_{a}$ and the shearing matrix $S_{s}$ (see \cite{Kuty})
\begin{align*}
M_{sa}=\left(\begin{array}{cccccc}
  a & \text{sgn}(a)\,a^{1/n}s_1 &  \text{sgn}(a)\, a^{1/n}s_2 & \text{sgn}(a)\, a^{1/n}s_3 & \cdots & \text{sgn}(a)\, a^{1/n}s_{n-1} \\
   0 & 0 & \text{sgn}(a)\, a^{1/n} & 0 & \cdots & \cdots \\
  \vdots & \vdots & \vdots& \vdots & \vdots & \vdots \\
   0 & 0 & 0 &0 & 0 & \text{sgn}(a)\, a^{1/n}
\end{array}\right).\tag{1.5}
\end{align*}

\pagestyle{myheadings}

\parindent=0mm \vspace{.0in}
The set $ \mathbb S={\mathbb R}\setminus \left\{0\right\}\times{\mathbb R^{n-1}}\times{\mathbb R^{n}}$ endowed with the operation
\begin{align*}
\big (a,s,t\big )\odot \big (a^{\prime},s^{\prime},t^{\prime}\big )=\big (aa^{\prime},s+a^{1-\frac{1}{n}} s^{\prime},t+S_{s}A_{a}t^{\prime}\big),\tag{1.6}
\end{align*}
forms a locally compact group, often called the {\it Shearlet group}. The left  Haar measures on $\mathbb S$ is given by  $d\eta={da\,ds\,dt}/a^{n+1}$ \cite{Dal2}. For every $\psi\in L^2(\mathbb R^n)$, we define
\begin{align*}
U(a,s,t)\psi(x)=\psi_{a,s,t}(x):=|\det M_{sa}|^{-1/2}\psi\big( M_{sa}^{-1}(x-t)\big).\tag{1.7}
\end{align*}
It is easy to verify that  $U: \mathbb S\to {\mathcal U}(L^{2}(\mathbb R^n))$ is a unitary mapping from the shearlet group $\mathbb S$ into the group of unitary operators ${\mathcal U}(L^{2}(\mathbb R^n))$ on $L^{2}(\mathbb R^n)$. In this framework, the continuous shearlet transform (1.1) takes the following form
\begin{align*}
{\mathcal {SH}}_{\psi}f(a,s,t)=\Big\langle{f,\psi_{a,s,t}}\Big\rangle=\Big\langle{f,U(a,s,t)\psi}\Big\rangle,\quad \text{for all}~f\in L^2(\mathbb R^n). \tag{1.8}
\end{align*}

\parindent=0mm \vspace{.0in}
The Heisenberg's uncertainty principle has played a fundamental role in the development and understanding of quantum mechanics,  signal processing and information theory \cite{Fol,HJ}. In quantum mechanics, this principle states that the position and the momentum of a particle cannot be both determined explicitly but only in a probabilistic sense with a certain degree of uncertainty. That is, increasing the knowledge of position, decreases the knowledge of momentum of the particle and vice-versa. The harmonic version of this principle says that a non-trivial function cannot be sharply localized in both time and frequency domains simultaneously \cite{CP,Bec}. With the development of time-frequency analysis, the study of  uncertainty principles have gained considerable attention and have been extended to a wide class of integral transforms including the short-time Fourier transform \cite{W}, fractional Fourier transform \cite{Guan},  wavelet transforms \cite{DM,Bat,Shah3}, fractional wavelet transform \cite{MA1}, linear canonical transforms \cite{Zhao} and special affine Fourier transforms \cite{Sun}.  The first study aimed to establish the uncertainty principles for the shearlet transforms was initiated by Dahlke et al.\cite{Dal3}, in which the authors have discussed various methods to minimize the uncertainty relations for the infinitesimal generators of the shearlet group. Later on, Su \cite{Su} derive some Heisenberg type uncertainty principles for the continuous shearlet transforms by adopting the strategy analogous to Wilcok \cite{W} and Cowling and Price \cite{CP}. Very recently, Nefzi et al.\cite{Nef} generalized the results of Su \cite{Su} for the multivariate shearlet transform and analyze the net concentration of these transforms on sets of finite measure using the machinery of projection operators. Recent results in this direction can be found in \cite{Shah4,BNT}.

\parindent=8mm \vspace{.1in}
To date, several generalizations, modifications and variations of the harmonic based uncertainty principles have appeared in the open literature, for instance,
the logarithmic uncertainty principles (Beckner-type uncertainty principles), entropy-based uncertainty relations, Benedick's uncertainty principles, Nazarov's uncertainty principles, local uncertainty principles and much more \cite{Bec,Ben,Vem,Naz,Jam,Kubo}. However, to the best of our knowledge, no such work has been explicitly carried out yet for the continuous shearlet transforms. It is therefore interesting and worthwhile to investigate these kinds of uncertainty principles for the continuous shearlet transforms in arbitrary space dimensions. The main objectives of this article are as follows:

\begin{itemize}
\item To obtain Pitt's inequality for the continuous shearlet transforms.
\item To establish Beckner's uncertainty principle for the continuous shearlet transforms.
\item To derive  Sobolev-type uncertainty inequalities for the continuous shearlet transforms.
\item To formulate Nazarov's uncertainty principle for the continuous shearlet transforms.
\item To obtain local uncertainty principles for the continuous shearlet transforms.
\end{itemize}

\parindent=0mm \vspace{.0in}
The rest of the article is structured as follows. In section 2, we establish an analogue of the well known Pitt's inequality for the continuous shearlet transforms in arbitrary space dimensions. In section 3, we derive the Beckner's uncertainty principle and obtain the corresponding Sobolev-type inequality for the continuous shearlet transforms. Sections 4 and 5 are respectively devoted to establishing the Nazarov's and local uncertainty principles for the shearlet transforms in arbitrary space dimensions. The conclusion is drawn in section 6.

\section{Pitt's Inequality for the Continuous Shearlet Transform}

The classical Pitt's inequality expresses a fundamental relationship between a sufficiently smooth function and the corresponding Fourier transform \cite{Bec}. For every $f\in \mathbb S(\mathbb R^n)\subseteq L^2(\mathbb R^n)$, the  inequality states that
\begin{align*}
\int_{\mathbb R^n}\left|\xi\right|^{-\lambda}\big|\mathscr F\big[f\big](\xi)\big|^2d\xi\le C_{\lambda}\int_{\mathbb R^n}\left|x\right|^{\lambda}\big|f(x)\big|^2dx,\quad 0\le\lambda<1\tag{2.1}
\end{align*}
where
\begin{align*}
C_{\lambda}=\pi^{\lambda}\left[\Gamma\left(\frac{n-\lambda}{4}\right)/ \Gamma\left(\frac{n+\lambda}{4}\right)\right]^2,\tag{2.2}
\end{align*}
and $\Gamma(\cdot)$ denotes the well known Euler's gamma function. Here, $\mathbb S(\mathbb R^n)$ denotes the Schwartz class in $L^2(\mathbb R^n)$ given by
\begin{align*}
\mathbb S\left(\mathbb R^n\right)=\left\{f\in C^{\infty}(\mathbb R^n): \sup_{t\in\mathbb R^n}\left|t^{\alpha}{\mathcal \partial}_{t}^{\beta}f(t)\right|<\infty\right\},\tag{2.3}
\end{align*}
where  $C^{\infty}(\mathbb R^n)$ is the class of smooth functions, $\alpha,\beta$ are any two non-negative integers, and ${\partial}_{t}$ denotes the usual partial differential operator.

\parindent=8mm \vspace{.1in}
The main objective of this section is to formulate an analogue of  Pitt's inequality (2.1) for the continuous shearlet transform in arbitrary space dimensions . Formally, we start our investigation with the following lemma.

\parindent=0mm \vspace{.1in}
{\bf Lemma 2.1.} {\it Let $\psi$ be an admissible shearlet, then for any $f\in{L^2(\mathbb R^{n})}$, we have }
\begin{align*}
{\mathscr F}\Big({\mathcal {SH}}_{\psi}f(a,s,t)\Big)\left(\xi\right)=\big|\det A_{a}\big|^{1/2}\hat{f}(\xi)\,\overline{\hat{\psi}\big( M_{sa}\xi\big)}. \tag{2.4}
\end{align*}

{\it Proof.} By virtue of Plancheral theorem for the classical Fourier transform, we obtain
\begin{align*}
{\mathcal {SH}}_{\psi}f(a,s,t)&=\int_{\mathbb R^{n}}f(x)~\overline{U(a,s,t)\,\psi(x)}\,dx\\
&=\int_{\mathbb R^{n}}\mathscr F\big[f\big](\xi)\,\overline{\mathscr F\Big[ U(a,s,t)\psi\Big]}(\xi)\,d{\xi}\\
&=\big|\det A_{a}\big|^{1/2}\int_{\mathbb R^{n}}\hat{f}(\xi)\left\{\int_{\mathbb R^{n}}\overline {\psi\big(D_{M_{sa}}(x-t)\big)\,e^{-2\pi i\xi \cdot x}\,dx}\right\}d\xi\\
&=\big|\det A_{a}\big|^{1/2}\int_{\mathbb R^{n}}\hat{f}(\xi)\left\{\int_{\mathbb R^{n}}\overline{\psi(z)\,e^{-2\pi i \xi\cdot(M_{sa}z+t)}\,dz}\right\}d\xi\\
&=\big|\det A_{a}\big|^{1/2}\int_{\mathbb R^{n}}\hat{f}(\xi)\,\overline{\hat{\psi}\big(M_{sa}\xi\big)\,e^{-2\pi i \xi\cdot t}}\,d\xi\\
&=\big|\det A_{a}\big|^{1/2}\int_{\mathbb R^{n}}\hat{f}(\xi)\,\overline{\hat{\psi}\big( M_{sa}\xi\big)}\,e^{2\pi i \xi\cdot t}\,d\xi\\
&=\big|\det A_{a}\big|^{1/2}{\mathscr F}^{-1}\Big[\hat{f}(\xi)\,\overline{\hat{\psi}\big( M_{sa}\xi\big)}\Big](\xi),
\end{align*}
which upon applying the Fourier transform yields
\begin{align*}
{\mathscr F}\Big({\mathcal {SH}}_{\psi}f(a,s,t)\Big)\left(\xi\right)=\big|\det A_{a}\big|^{1/2}\hat{f}(\xi)\,\overline{\hat{\psi}\big( M_{sa}\xi\big)}.
\end{align*}

This completes the proof of Lemma 2.1. \quad \fbox

\parindent=8mm \vspace{.1in}

 In our next lemma, we shall establish  the Moyal's principle for the continuous shearlet transform (1.1) in arbitrary space dimensions, which will  be employed in the subsequent sections to obtain certain uncertainty inequalities.

\parindent=0mm \vspace{.1in}

{\bf Lemma 2.2.} {\it Let $\big[{\mathcal {SH}}_{\psi}f\big](a,s,t)$ and $\big[{\mathcal {SH}}_{\psi}g\big](a,s,t)$ be the shearlet transforms for a given pair of square integrable functions $f$ and $g$. Then,  the following identity holds:}
\begin{align*}
\int_{\mathbb S}\Big({\mathcal {SH}}_{\psi}f(a,s,t)\Big)\overline{\Big({\mathcal {SH}}_{\psi}g(a,s,t)\Big)}\,d\eta=C_{\psi}\, \big\langle f,g \big\rangle, \tag{2.5}
\end{align*}
{\it where $C_{\psi}$ is the admissability condition of the shearlet $\psi\in L^2(\mathbb R^n)$ given by}
\begin{align*}
C_{\psi}=\int_{\mathbb R^{n-1}}\int_{{\mathbb R \setminus \left\{0\right\}}}\dfrac{\big|{\hat{\psi}( M_{sa}\xi)}\big|^2}{a^{\frac{n^2-n+1}{n}}}\,da\,ds<\infty.\tag{2.6}
\end{align*}
{\it Proof.} Using the unitary representation of the continuous   shearlet transform (1.8), we have
\begin{align*}
{\mathcal {SH}}_{\psi}f(a,s,t)&={\Big\langle f, U(a,s,t)\psi \Big\rangle}\\
&={\Big\langle \mathscr F\big[f\big](\xi),\mathscr F \big[U(a,s,t)\psi\big](\xi)\Big\rangle}\\
&=\int_{\mathbb R^{n}}\mathscr F\big[f\big](\xi)~\overline{\mathscr F \Big[U(a,s,t)\psi\Big]}(\xi)\,d\xi\\
&=\big|\det A_{a}\big|^{1/2}\int_{\mathbb R^{n}}\hat{f}(\xi)\,\overline{\hat{\psi}\big(M_{sa}\xi \big)}\,e^{2\pi i \xi\cdot t}\,d\xi.\tag{2.7}
\end{align*}
Similarly, we have
\begin{align*}
{\mathcal {SH}}_{\psi}g(a,s,t)=\big|\det A_{a}\big|^{1/2}\int_{\mathbb R^{n}}\hat{g}(\sigma)\,\overline{\hat{\psi}\big(M_{sa}\sigma \big)}\,e^{2\pi i \sigma \cdot t}\,d\sigma.\tag{2.8}
\end{align*}
An implication of the well-known Fubini theorem yields
\begin{align*}
&\int_{\mathbb S}\Big({\mathcal {SH}}_{\psi}f(a,s,t)\Big)\overline{\Big({\mathcal {SH}}_{\psi}g(a,s,t)\Big)}\,d\eta\\
&=\int_{\mathbb R^{n}}\int_{\mathbb R^{n-1}}\int_{\mathbb R \setminus \left\{0\right\}}\Big({\mathcal {SH}}_{\psi}f(a,s,t)\Big)\overline{\Big({\mathcal {SH}}_{\psi}g(a,s,t)\Big)}\,\dfrac{\,da\,ds\,dt}{a^{n+1}}\\
&=\int_{\mathbb R^{n}}\int_{\mathbb R^{n-1}}\int_{\mathbb R \setminus \left\{0\right\}}\left\{\int_{\mathbb R^{n}}\int_{\mathbb R^{n}}\big|\det A_{a}\big|\hat{f}(\xi)\,\overline{\hat{\psi}( M_{sa}\xi)}\,e^{2\pi i \xi\cdot t}\,  \overline{\hat{g}(\sigma)}\,\,\hat{\psi}\big( M_{sa}\sigma\big)\,e^{-2\pi i \sigma \cdot t} d\xi\, d\sigma\right\}\dfrac{\,da\,ds\,dt}{a^{n+1}}\\
&=\int_{\mathbb R^{n-1}}\int_{\mathbb R \setminus \left\{0\right\}}\dfrac{\,da\,ds}{a^{\frac{n^2-n+1}{n}}}\int_{\mathbb R^{n}}\int_{\mathbb R^{n}}\hat{f}(\xi)\,\overline{\hat{\psi}( M_{sa}\xi)}\,  \,\overline{\hat{g}(\sigma)}\,\hat{\psi}\big( M_{sa}\sigma\big)\,\left\{\int_{\mathbb R^{n}}e^{2\pi i(\xi-\sigma)\cdot\, t}\, dt\right\}d\xi\, d\sigma\\
&=\int_{\mathbb R^{n-1}}\int_{\mathbb R \setminus \left\{0\right\}}\dfrac{\,da\,ds}{a^{\frac{n^2-n+1}{n}}}\int_{\mathbb R^{n}}\int_{\mathbb R^{n}}\hat{f}(\xi)\,\overline{\hat{\psi}( M_{sa}\xi)}\,\,\overline{\hat{g}(\sigma)}\,\hat{\psi}\big( M_{sa}\sigma\big)\,\,\delta(\sigma-\xi)\,d\xi\, d\sigma\\
&=\int_{\mathbb R^{n-1}}\int_{\mathbb R \setminus \left\{0\right\}}\dfrac{\,da\,ds}{a^{\frac{n^2-n+1}{n}}}\int_{\mathbb R^{n}}\hat{f}(\xi)\,\overline{\hat{g}(\xi)}\,\,\left|{\hat{\psi}\big( M_{sa}\xi\big)}\right|^2\,d\xi\\
&=\int_{\mathbb R^{n}}\hat{f}(\xi)\,\overline{\hat{g}(\xi)}\,\left\{\int_{\mathbb R^{n-1}}\int_{\mathbb R \setminus \left\{0\right\}}\dfrac{\left|{\hat{\psi}\big( M_{sa}\xi\big)}\right|^2}{a^{\frac{n^2-n+1}{n}}}\,da\,ds\right\}\,d\xi\\
&=C_{\psi}\left\langle \hat{f},\hat{g}\right\rangle\\
&=C_{\psi}\,\big\langle f, g\big\rangle.
\end{align*}

This completes the proof of Lemma 2.2. \quad \fbox

\parindent=0mm \vspace{.1in}

{\it Remarks.} (i). For $f=g$, equation (2.5) yields the following energy preserving relation
\begin{align*}
&\int_{\mathbb S}\Big|{\mathcal {SH}}_{\psi}f(a,s,t)\Big|^{2}{d\eta}=C_{\psi}\big\|f\big\|^{2}_{2}.\tag{2.9}
\end{align*}
(ii). Equation (2.9) demonstrates  that the continuous shearlet transform (1.1) is a bounded linear operator from $L^2(\mathbb R^n)$ to $L^2({\mathbb R \setminus \left\{0\right\}}\times\mathbb R^{n-1}\times\mathbb R^n)$.

\parindent=0mm \vspace{.1in}
(iii). For $C_{\psi}=1$,  the continuous shearlet transform (1.1) becomes an isometry from the space of signals $L^2(\mathbb R^n)$ to the space of transforms $L^2({\mathbb R \setminus \left\{0\right\}}\times\mathbb R^{n-1}\times\mathbb R^n).$

\parindent=8mm \vspace{.1in}
We are now in a position to establish the Pitt's inequality  for the continuous shearlet transforms in arbitrary space dimensions.

\parindent=0mm \vspace{.1in}

{\bf Theorem 2.3.} {\it For any arbitrary $f\in \mathcal S(\mathbb R^n)\subseteq L^2(\mathbb R^n)$, the Pitt's inequality for the continuous shearlet transform (1.1) is given by:}
\begin{align*}
C_{\psi}\int_{\mathbb R^{n}}\left|\xi\right|^{-\lambda}\left|\hat{f}(\xi)\right|^2\,d\xi\le C_{\lambda}\int_{\mathbb S}\left|t\right|^{\lambda}\Big|\mathcal {SH}_{\psi}f(a,s,t)\Big|^2 d\eta.\tag{2.10}
\end{align*}
{\it where $C_{\psi}$ is the admissability condition given by (2.6).}

\parindent=0mm \vspace{.1in}
{\it Proof.} As a consequence of the  inequality (2.1), we can write
\begin{align*}
\int_{\mathbb R^n}\left|\xi\right|^{-\lambda}\left|\mathscr F\Big[\mathcal {SH}_{\psi}f(a,s,t)\Big](\xi)\right|^2d\xi\le C_{\lambda}\int_{\mathbb R^n}\left|t\right|^{\lambda}\Big|\mathcal {SH}_{\psi}f(a,s,t)\Big|^2dt,\tag{2.11}
\end{align*}
which upon integration  with respect to the Haar measure $d\eta=dsda/a^{n+1}$ yields
\begin{align*}
\int_{\mathbb R^{n}}\int_{\mathbb R^{n-1}}\int_{{\mathbb R \setminus \left\{0\right\}}}\left|\xi\right|^{-\lambda}\Big|\mathscr F\Big[\mathcal {SH}_{\psi}f(a,s,t)\Big](\xi)\Big|^2\dfrac{d\xi\,ds\,da}{a^{n+1}}\le C_{\lambda}\int_{\mathbb S}\left|t\right|^{\lambda}\Big|\mathcal {SH}_{\psi}f(a,s,t)\Big|^2 d\eta.\tag{2.12}
\end{align*}
Invoking Lemma 2.1, we can express the inequality (2.12) in the following manner:
\begin{align*}
\int_{\mathbb R^{n}}\int_{\mathbb R^{n-1}}\int_{{\mathbb R \setminus \left\{0\right\}}}\left|\xi\right|^{-\lambda}\left|\hat{f}(\xi)~\overline{\hat{\psi}\big( M_{sa}\xi\big)}\right|^2\,\dfrac{d\xi\,ds\,da}{a^{\frac{n^2-n+1}{n}}}\le C_{\lambda}\int_{\mathbb S}\left|t\right|^{\lambda}\Big|\mathcal {SH}_{\psi}f(a,s,t)\Big|^2 d\eta.\tag{2.13}
\end{align*}
Equivalently, we have
\begin{align*}
\int_{\mathbb R^{n}}\left|\xi\right|^{-\lambda}\left|\hat{f}(\xi)\right|^2\left\{\int_{\mathbb R^{n-1}}\int_{{\mathbb R \setminus \left\{0\right\}}}\dfrac{\big|{\hat{\psi}( M_{sa}\xi)}\big|^2}{a^{\frac{n^2-n+1}{n}}}\,da\,ds \right\}\,d\xi\le C_{\lambda}\int_{\mathbb S}\left|t\right|^{\lambda}\Big|\mathcal {SH}_{\psi}f(a,s,t)\Big|^2 d\eta.\tag{2.14}
\end{align*}
Since $\psi$ is an admissible shearlet, therefore inequality (2.14) becomes
\begin{align*}
C_{\psi}\int_{\mathbb R^{n}}\left|\xi\right|^{-\lambda}\left|\hat{f}(\xi)\right|^2\,d\xi\le C_{\lambda}\int_{\mathbb S}\left|t\right|^{\lambda}\Big|\mathcal {SH}_{\psi}f(a,s,t)\Big|^2 d\eta,\tag{2.15}
\end{align*}
which establishes the Pitt's inequality for the continuous shearlet transform in arbitrary space dimensions. \quad \fbox

\parindent=0mm \vspace{.1in}

{\it Remark:}  For $\lambda=0$, equality  holds in (2.10), which is in consonance with the classical Pitt's inequality.

\section{Beckner-type Inequalities for the Continuous Shearlet Transforms}

\parindent=0mm \vspace{.1in}
The classical Beckner's inequality \cite{Bec} is given by
\begin{align*}
\int_{\mathbb R^{n}}{\ln|t|} ~{\big|f(t)\big|^{2}}\,dt+\int_{\mathbb R^{n}}{\ln{|\xi|}\left|\hat{f}(\xi)\right|^2}\,d{\xi}\geq \left(\dfrac{{\Gamma^{\prime}}(1/2)}{\Gamma(1/2)}-\ln{\pi}\right)\int_{\mathbb R^{n}} \big|f(t)\big|^{2}\,dt\tag{3.1}
\end{align*}
for all $f\in L^{2}(\mathbb R^{n})$, for which the quantity on left is defined, where $t\in{\mathbb R^{2}}$, and $\Gamma(t)$ is the  gamma function. This inequality is related to the classical Heisenberg's uncertainty principle and for that reason it is often referred as the logarithmic uncertainty principle. Considerable attention has been paid to this inequality for its various generalizations, improvements, analogues,  and their applications in science and engineering (see \cite{Fol,HJ,Guan,MA1,Shah4}).

\parindent=0mm \vspace{.1in}

{\bf Theorem 3.1.} {\it Let $\big[{\mathcal {SH}}_{\psi}f\big](a,s,t)$  be the shearlet transform of any  arbitrary function $f\in \mathbb S(\mathbb R^n)$, the following inequality holds:}
\begin{align*}
\int_{\mathbb S}{\ln|t|} \,{\Big|{\mathcal {SH}}_{\psi}f(a,s,t)\Big|^{2}}\,d\eta+{C_{\psi}}\int_{\mathbb R^{n}}{\ln{|\xi|} \left|\hat{f}(\xi)\right|^2}\,d{\xi}\geq {C_{\psi}}\left[\dfrac{{\Gamma^\prime}(n/4)}{\Gamma(n/4)}-\ln{\pi}\right]\big\|f\big\|_{2}^{2},\tag{3.2}
\end{align*}
{\it where $C_{\psi}$ is given by (2.6).}

\parindent=0mm \vspace{.2in}
{\it Proof.} For  every $0\le \lambda <1$, we define
\begin{align*}
P\left({\lambda}\right)=C_{\psi}\int_{\mathbb R^{n}}\left|\xi\right|^{-\lambda}\left|\hat{f}(\xi)\right|^2\,d\xi- C_{\lambda}\int_{\mathbb S}\left|t\right|^{\lambda}\Big|\mathcal {SH}_{\psi}f(a,s,t)\Big|^2 d\eta.\tag{3.3}
\end{align*}
On differentiating (3.3) with respect to $\lambda$, we obtain
\begin{align*}
P^{\prime}\left({\lambda}\right)=-C_{\psi}\int_{\mathbb R^{n}}\left|\xi\right|^{-\lambda}\ln\big|\xi\big|\left|\hat{f}(\xi)\right|^2\,d\xi- C_{\lambda}\int_{\mathbb S}\left|t\right|^{\lambda}\ln\big|t\big|\Big|\mathcal {SH}_{\psi}f(a,s,t)\Big|^2 d\eta\\
-C^{\,\prime}_{\lambda}\int_{\mathbb S}\left|t\right|^{\lambda}\Big|\mathcal {SH}_{\psi}f(a,s,t)\Big|^2 d\eta.\tag{3.4}
\end{align*}
where
\begin{align*}
C^{\,\prime}_{\lambda}&=-\dfrac{\pi^{\lambda}}{2}\left\{\dfrac{{\Gamma}^2\left(\dfrac{n+\lambda}{4}\right)\Gamma\left(\dfrac{n-\lambda}{4}\right)
{\Gamma}^{\prime}\left(\dfrac{n-\lambda}{4}\right)+{\Gamma}^2\left(\dfrac{n-\lambda}{4}\right)\Gamma\left(\dfrac{n+\lambda}{4}\right)
{\Gamma}^{\prime}\left(\dfrac{n+\lambda}{4}\right)}{{\Gamma}^2\left(\dfrac{n+\lambda}{4}\right)}\right\}\\
&\qquad\qquad\qquad\qquad\qquad\qquad\qquad\qquad+\pi^{\lambda}\ln \pi
\left\{{\Gamma}^2\left(\dfrac{n-\lambda}{4}\right)\big /{\Gamma}^2\left(\dfrac{n+\lambda}{4}\right)\right\}.\tag{3.5}
\end{align*}
For $\lambda=0$, equation (3.5) yields
\begin{align*}
C^{\,\prime}_{0}=\left[\ln \pi-\dfrac{\Gamma^{\prime}(n/4)}{\Gamma (n/4)}\right].\tag{3.6}
\end{align*}
By virtue of Pitt's inequality (2.10) for the shearlet transforms, it follows that $P(\lambda)\le0$, for all  $\lambda\in [0,1)$ and
\begin{align*}
P(0)&=C_{\psi}\int_{\mathbb R^{n}}\left|\hat{f}(\xi)\right|^2\,d\xi- C_{0}\int_{\mathbb S}\Big|\mathcal {SH}_{\psi}f(a,s,t)\Big|^2 d\eta
=C_{\psi}\left\|\hat{f}\right\|_2^2-C_{\psi}\big\|f\big\|_2^2=0.\tag{3.7}
\end{align*}
Therefore, for any $h>0$, we observe that $P^{\,\prime}\left(0+h\right)\le0$, whenever $h\rightarrow0$; that is,
\begin{align*}
&-C_{\psi}\int_{\mathbb R^{n}}\ln\big|\xi\big|\left|\hat{f}(\xi)\right|^2d\xi- C_{0}\int_{\mathbb S}\ln\big|t\big|\Big|\mathcal {SH}_{\psi}f(a,s,t)\Big|^2 d\eta
-C^{\,\prime}_{0}\int_{\mathbb S}\Big|\mathcal {SH}_{\psi}f(a,s,t)\Big|^2 d\eta\le 0.\tag{3.8}
\end{align*}
Applying the energy preserving relation (2.9) and the obtained estimate (3.6) of $C^{\,\prime}_{0}$, we obtain
\begin{align*}
-C_{\psi}\int_{\mathbb R^{n}}\ln\big|\xi\big|\left|\hat{f}(\xi)\right|^2\,d\xi- \int_{\mathbb S}\ln\big|t\big|\Big|\mathcal {SH}_{\psi}f(a,s,t)\Big|^2 d\eta
-\left[\ln \pi-\dfrac{\Gamma^{\prime}(1/2)}{\Gamma (1/2)}\right]C_{\psi}\,\big\|f\big\|_2^2\le0,
\end{align*}
or equivalently,
\begin{align*}
\int_{\mathbb S}\ln\big|t\big|\Big|\mathcal {SH}_{\psi}f(a,s,t)\Big|^2 d\eta+
C_{\psi}\int_{\mathbb R^{n}}\ln\big|\xi\big|\left|\hat{f}(\xi)\right|^2 d\xi\ge
\left[\dfrac{\Gamma^{\prime}(n/4)}{\Gamma (n/4)}-\ln \pi\right]C_{\psi}\,\big\|f\big\|_2^2.\tag{3.9}
\end{align*}
Inequality (3.9) is the desired Beckner's uncertainty principle for the continuous shearlet transform in arbitrary space dimensions.\quad\fbox

\parindent=8mm \vspace{.1in}

We now present an alternate proof of Theorem 3.1. The strategy of the proof is different and is obtained directly from the  classical Beckner's inequality (3.1).

\parindent=0mm \vspace{.1in}
{\it Second Proof of the Theorem 3.1.} We shall identify ${\mathcal {SH}}_{\psi}f(a,s,t)$ as a function of the translation parameter $t$ and then replace $f \in \mathbb S(\mathbb R^n)$ in (3.1) with  ${\mathcal {SH}}_{\psi}f(a,s,t)$, so that
\begin{align*}
\int_{\mathbb R^{n}}{\ln|t|} \,{\Big|{\mathcal {SH}}_{\psi}f(a,s,t)\Big|^{2}}\,dt&+\int_{\mathbb R^{n}}{\ln{|\xi|}\Big|{\mathscr F}\Big[{\mathcal {SH}}_{\psi}f(a,s,t)\Big](\xi)\Big|^2}\,d{\xi}\\
&\ge\left(\dfrac{{\Gamma^{\prime}}(n/4)}{\Gamma(n/4)}-\ln{\pi}\right)\int_{\mathbb R^{n}} \Big|{\mathcal {SH}}_{\psi}f(a,s,t)\Big|^{2}\,dt.\tag{3.10}
\end{align*}
Integrating (3.10) with respect to the measure $d\eta=dads/a^{n+1}$, we obtain
\begin{align*}
\int_{\mathbb S}  {\ln|t|} \,{\Big|{\mathcal {SH}}_{\psi}f(a,s,t)\Big|^{2}} d\eta+\int_{{\mathbb R \setminus \left\{0\right\}}}\int_{\mathbb R}\int_{\mathbb R^{n}}{\ln{|\xi|}\, \Big|{\mathscr F}\Big[{\mathcal {SH}}_{\psi}f(a,s,t)\Big](\xi)\Big|^2}\,\dfrac{da\,ds\,d{\xi}}{a^{n+1}}\\
\geq  \left(\dfrac{{\Gamma^{\prime}}(n/4)}{\Gamma(n/4)}-\ln{\pi}\right)\int_{\mathbb S}\Big|{\mathcal {SH}}_{\psi}f(a,s,t)\Big|^{2} d\eta.\tag{3.11}
\end{align*}
Using equation (2.9), we have
\begin{align*}
\int_{\mathbb S}\,\ln|t| \,{\Big|{\mathcal {SH}}_{\psi}f(a,s,t)\Big|^{2}}\,d\eta+\int_{\mathbb S}\,{\ln{|\xi|}\, \Big|{\mathscr F}\Big[{\mathcal {SH}}_{\psi}f(a,s,t)\Big](\xi)\Big|^2}\,d\eta\geq  \left(\dfrac{{\Gamma^{\prime}}(n/4)}{\Gamma(n/4)}-\ln{\pi}\right)C_{\psi}\big\|f\big\|_2^2.\tag{3.12}
\end{align*}
We shall now simplify the second integral  of (3.12) as
\begin{align*}
&\int_{\mathbb S}\ln|\xi| \Big|{\mathscr F}\Big[{\mathcal {SH}}_{\psi}f(a,s,t)\Big]\left(\xi\right)\Big|^2 d\eta \\
&\quad=\int_{\mathbb R^{n-1}}\int_{{\mathbb R \setminus \left\{0\right\}}}\int_{\mathbb R^{n}}\ln|\xi|\Big|{\mathscr F}\Big[{\mathcal {SH}}_{\psi}f(a,s,t)\Big]\left(\xi\right)\Big|^2 \dfrac{d{\xi}\,da\,ds}{a^{n+1}}\\
&\quad=\int_{\mathbb R^{n-1}}\int_{{\mathbb R \setminus \left\{0\right\}}}\int_{\mathbb R^{n}}\ln|\xi|\Big[{\mathscr F}\Big[{\mathcal {SH}}_{\psi}f(a,s,t)\Big]\left(\xi\right)\Big]\overline{\Big[{\mathscr F}\Big[{\mathcal {SH}}_{\psi}f(a,s,t)\Big]\left(\xi\right)\Big]}\,\dfrac{d{\xi}\,da\,ds}{a^{n+1}}\\
&\quad=\int_{\mathbb R^{n-1}}\int_{{\mathbb R \setminus \left\{0\right\}}}\int_{\mathbb R^{n}}\ln\left|\xi\right| \big|\det A_{a}\big|\,\hat{f}(\xi)\,\overline{\hat{\psi}\big( M_{sa}\xi\big)}\,\overline{\hat{f}(\xi)}\,\hat{\psi}\big(\xi M_{sa}\big)\,\dfrac{d{\xi}\,da\,ds}{a^{n+1}}\\
&\quad=\int_{\mathbb R^{n-1}}\int_{{\mathbb R \setminus \left\{0\right\}}}\int_{\mathbb R^{n}}\big|\det A_{a}\big|\,\ln\left|\xi\right|\, \left|\hat{\psi}( M_{sa}\xi)\right|^{2}\,\left|\hat{f}(\xi)\right|^{2}\dfrac{d{\xi}\,da\,ds}{a^{n+1}}\\
&\quad=\int_{\mathbb R^{n-1}}\int_{{\mathbb R \setminus \left\{0\right\}}}\int_{\mathbb R^{n}}\big|\det A_{a}\big|\ln\left|\xi\right| \left|\hat{\psi}( M_{sa}\xi)\right|^{2}\left|\hat{f}(\xi)\right|^{2}\dfrac{d{\xi}\,da\,db}{a^{n+1}}\\
&\quad=\int_{\mathbb R^{n}}{\ln{\left|\xi\right|} \left|\hat{f}(\xi)\right|^2}\Bigg\{\int_{\mathbb R^{n-1}}\int_{{\mathbb R \setminus \left\{0\right\}}}\dfrac{\big|\hat{\psi}(M_{sa}\xi)\big|^{2}}{a^{\frac{n^2-n+1}{n}}}\,da\,ds\Bigg\}d{\xi}\\
&\quad=C_{\psi}\int_{\mathbb R^{n}}{\ln{\left|\xi\right|} \left|\hat{f}(\xi)\right|^2}d{\xi}.\tag{3.13}
\end{align*}
Plugging the estimate (3.13) in (3.12) gives the desired inequality for the continuous shearlet transforms as
\begin{align*}
\int_{\mathbb S}{\ln|t|} \,{\Big|{\mathcal {SH}}_{\psi}f(a,s,t)\Big|^{2}}\,d\eta + C_{\psi}\int_{\mathbb R^{n}}\ln|\xi| \left|\hat{f}(\xi)\right|^{2} d{\xi}\geq \left(\dfrac{{\Gamma^{\prime}}(n/4)}{\Gamma(n/4)}-\ln{\pi}\right)C_{\psi}\big\|f\big\|^{2}_{2}.
\end{align*}
This completes the second proof of Theorem 3.1. \quad \fbox

\parindent=0mm \vspace{.1in}

{\it{Deduction:}} Using Jensen's inequality in (3.2), we  obtain  an analogue of the classical Heisenberg's uncertainty inequality  for the  continuous shearlet transforms  as
\begin{align*}
&\ln \left\{ \int_{\mathbb S} |t|^2 \,\Big|{\mathcal {SH}}_{\psi}f(a,s,t)\Big|^2 d\eta \, C_{\psi} \int_{\mathbb R^n} |\xi|^2 \left|\hat f(\xi)\right|^2 d\xi \right\}^{1/2}\\
&= \ln \left\{ \int_{\mathbb S} |t|^2 \,\Big|{\mathcal {SH}}_{\psi}f(a,s,t)\Big|^2 d\eta\right\}^{1/2}+\ln \left(C_{\psi}\right)^{1/2}+\ln \left\{ \int_{\mathbb R^n} |\xi|^2 \left|\hat f(\xi)\right|^2 d\xi \right\}^{1/2}\\
&\ge \int_{\mathbb S}{\ln|t|} \,{\Big|{\mathcal {SH}}_{\psi}f(a,s,t)\Big|^{2}}\,d\eta +\ln \left(C_{\psi}\right)^{1/2}+\int_{\mathbb R^{n}}\ln|\xi| \left|\hat{f}(\xi)\right|^{2} d{\xi}\\
&\ge \left(\dfrac{{\Gamma^{\prime}}(n/4)}{\Gamma(n/4)}-\ln{\pi}\right)C_{\psi}\big\|f\big\|^{2}_{2}+\ln \left(C_{\psi}\right)^{1/2},
\end{align*}
which upon simplification with $C_\psi=1$ yields
\begin{align*}
\left\{\int_{\mathbb S} |t|^2 \,\Big|{\mathcal {SH}}_{\psi}f(a,s,t)\Big|^2 d\eta\right\}^{1/2}\left\{  \int_{\mathbb R^n} |\xi|^2 \left|\hat f(\xi)\right|^2 d\xi \right\}^{1/2}\ge \exp\left\{\dfrac{ -2\sqrt{\pi}\ln 2}{\sqrt {\pi}}-\ln \pi\right\}\big\|f\big\|^{2}_{2}= \dfrac{\big\|f\big\|^{2}_{2}}{4\pi}.
\end{align*}

\parindent=0mm \vspace{.1in}
The remaining part of this Section is devoted to establish the Sobolev-type uncertainty inequality for the continuous shearlet transform in arbitrary space dimensions. This inequality is employed in the Section 5 to obtain a local-type uncertainty principle for the continuous shearlet transform (1.1). To facilitate our intention, we start with the following definitions:

\parindent=0mm \vspace{.1in}
{\bf Definition 3.2.} The Sobolev space on $\mathbb R^n$ is  defined by
\begin{align*}
\mathbb H\left(\mathbb R^n\right)=\Big\{f\in L^2(\mathbb R^n): \nabla f\in L^2(\mathbb R^n)\Big\},\tag{3.14}
\end{align*}
where $\nabla$ denotes the differential operator given by $\nabla=\left(\dfrac{\partial}{\partial x_1},\dfrac{\partial}{\partial x_2},\dots,\dfrac{\partial}{\partial x_n}\right)$.

\parindent=0mm \vspace{.1in}
{\bf Definition 3.3.} For $1\le p<\infty$ and $b>0$,  the weighted Lebesgue space on $\mathbb R^n$ is defined by
\begin{align*}
\mathbb W_{b}^{p}\left(\mathbb R^n\right)=\Big\{f\in L^p(\mathbb R^n): \langle t\rangle^{b}f\in L^p(\mathbb R^n)\Big\},\tag{3.15}
\end{align*}
where $ \langle t\rangle$ is the weight function given by  $\langle t\rangle=\big(1+\left|t\right|^2\big)^{1/2}$, $t\in\mathbb R^n$.

\parindent=8mm \vspace{.1in}
The logarithmic Sobolev inequality is related to the class of functions $\mathbb H\left(\mathbb R^n\right)$ states that for any non-trivial function $f\in\mathbb H\left(\mathbb R^n\right)$ \cite{PD},
\begin{align*}
\int_{\mathbb R^n}\big|f(t)\big|^2\ln\left(\dfrac{\left|f(t)\right|^2}{\left\|f\right\|_2^2}\right)\,dt\le
\dfrac{n}{2}\ln\left(\dfrac{2}{n\pi e \left\|f\right\|_2^2}\int_{\mathbb R^n}\big|\nabla f(t)\big|^2dt\right).\tag{3.16}
\end{align*}
Inequality (3.18) is often referred as Gross's inequality \cite{PD,Gen}. On the other hand, Beckner \cite{Bec} proved another version of logarithmic Sobolev inequality for extremal functions which offers  better estimate than Gross's inequality (3.16) given by
\begin{align*}
\int_{\mathbb R^n}\big|f(t)\big|^2\ln\left(\dfrac{\left|f(t)\right|^2}{\left\|f\right\|_2^2}\right)\,dt\le
\dfrac{n}{2}\int_{\mathbb R^n}\left|\hat{f}(\xi)\right|^2\ln \left(B_n\left|\xi\right|^2\right)\,d\xi- n\big\|f\big\|^2_{2} \left(\dfrac{\Gamma^{\prime}(n/2)}{\Gamma(n/2)}\right),\tag{3.17}
\end{align*}
where $B_n=\frac{1}{4\pi}\left(\frac{\Gamma (n)}{\Gamma(n/2)}\right)^{2/n}$.

\parindent=8mm \vspace{.1in}
Very recently, Kubo et al.\cite{Kubo} obtained a logarithmic Sobolev-type inequality for the weighted Lebesgue spaces $\mathbb W_{b}^{p}\left(\mathbb R^n\right)$ and pointed out that the obtained inequality has a dual relation with the Beckner's inequality (3.1).   For any non-trivial function $f\in \mathbb W_{b}^{1}\big({\mathbb R^n}\big),$ the inequality states that
\begin{align*}
-\int_{\mathbb R^n}\big|f(t)\big|\ln\left\{\dfrac{\left|f(t)\right|}{\left\|f\right\|_1}\right\}\,dt\le
n\int_{\mathbb R^n}\big|f(t)\big| \ln \left\{C_{n,b}\left(1+\left|t\right|^{b}\right)\right\}\,dt,\tag{3.18}
\end{align*}
where
\begin{align*}
C_{n,b}=\left\{\dfrac{2\,\pi^{n/2}\Gamma(n/b)\,\Gamma(n/b^{\prime})}{b\, \Gamma(n)\,\Gamma(n/2)}\right\}^{1/n}\quad {\text {and}}\quad \dfrac{1}{b}+\dfrac{1}{b^{\,\prime}}=1.\qquad\qquad\tag{3.19}
\end{align*}

\parindent=0mm \vspace{.0in}
Furthermore, the duality has been shown in the following sense:
\begin{align*}
\int_{\mathbb R^{n}}\big|f(t)\big|^2\ln\left(\dfrac{1+\left|t\right|^2}{2}\right)\,dt+\int_{\mathbb R^{n}}\left|\hat {f}(\xi)\right|^2\ln\big|\xi\big|\,d{\xi}\ge\left(\dfrac{{\Gamma^{\prime}}(n/2)}{\Gamma(n/2)}\right)\int_{\mathbb R^{n}}\big|f(t)\big|^2 dt.\tag{3.20}
\end{align*}

\parindent=0mm \vspace{.0in}

The following theorem is the main result of this subsection which establishes an analogue of the Sobolev-type uncertainty inequality (3.20) for the continuous shearlet transforms in arbitrary space dimensions.

\parindent=0mm \vspace{.1in}

{\bf Theorem 3.4.} {\it If $\big[{\mathcal {SH}}_{\psi}f\big](a,s,t)$  is the shearlet transform of any  arbitrary function $f \in \mathbb H(\mathbb R^n)\cap \mathbb W_1^{1}(\mathbb R^n)$, then the following Sobolev-type uncertainty inequality holds:}
\begin{align*}
\int_{\mathbb S} \Big|{\mathcal {SH}}_{\psi}f(a,s,t)\Big|^2\ln\left(\dfrac{1+\left|t\right|^2}{2}\right)\,d\eta+C_{\psi}\int_{\mathbb R^{n}}{\ln{\left|\xi\right|} \left|\hat{f}(\xi)\right|^2}d{\xi}
\ge\left(\dfrac{{\Gamma^{\prime}}(n/2)}{\Gamma(n/2)}\right)C_{\psi}\big\|f\big\|^{2}_{2},\tag{3.21}
\end{align*}
{\it whenever the L.H.S of (3.21) is defined.}

\parindent=0mm \vspace{.1in}
{\it Proof.} As a consequence of inequality (3.20), we have
\begin{align*}
\int_{\mathbb R^{n}}\Big|{\mathcal {SH}}_{\psi}f(a,s,t)\Big|^2\ln\left(\dfrac{1+\left|t\right|^2}{2}\right)\,dt+\int_{\mathbb R^{n}}\Big|\mathscr F \Big[{\mathcal {SH}}_{\psi}f(a,s,t)\Big](\xi)\Big|^2\ln\big|\xi\big|\,d{\xi}\\
\qquad\qquad\ge\left(\dfrac{{\Gamma^{\prime}}(n/2)}{\Gamma(n/2)}\right)\int_{\mathbb R^{n}}\Big|{\mathcal {SH}}_{\psi}f(a,s,t)\Big|^2dt,
\end{align*}
which upon integration yields
\begin{align*}
\int_{\mathbb S} \Big|{\mathcal {SH}}_{\psi}f(a,s,t)\Big|^2\ln\left(\dfrac{1+\left|t\right|^2}{2}\right)\,d\eta+\int_{\mathbb S}\ln\big|\xi\big|\, \Big|\mathscr F \Big[{\mathcal {SH}}_{\psi}f(a,s,t)\Big](\xi)\Big|^2\,d\eta\\
\qquad\qquad\ge\left(\dfrac{{\Gamma^{\prime}}(n/2)}{\Gamma(n/2)}\right)\int_{\mathbb S}\Big|{\mathcal {SH}}_{\psi}f(a,s,t)\Big|^2d\eta.\tag{3.22}
\end{align*}
Using the estimate  (3.13) for the second integral on the L.H.S of (3.22) and invoking (2.9), we obtain
\begin{align*}
\int_{\mathbb S} \Big|{\mathcal {SH}}_{\psi}f(a,s,t)\Big|^2\ln\left(\dfrac{1+\left|t\right|^2}{2}\right)\,d\eta+C_{\psi}\int_{\mathbb R^{n}}{\ln{\left|\xi\right|} \left|\hat{f}(\xi)\right|^2}d{\xi}
\ge\left(\dfrac{{\Gamma^{\prime}}(n/2)}{\Gamma(n/2)}\right)C_{\psi}\big\|f\big\|^{2}_{2},\tag{3.23}
\end{align*}
where $C_{\psi}$ is given by (2.6). This completes the proof of Theorem 3.4.\quad \fbox

\section{Nazarov-type Inequality for the Shearlet Transforms}

As is well known, the classical Heisenberg's uncertainty principle measures the localization in terms of the dispersions of the respective functions. Considering an alternate criterion of localization; that is, the smallness of the support, Nazarov \cite{Naz,Jam}  proposed an uncertainty principle which is concerned with the query;  what happens if a non-zero function and its Fourier transform are small outside a compact set? The Nazarov's uncertainty principle in the classical Fourier domain states that if $E_1$ and $E_2$ are two subsets of $\mathbb R^n$ with finite measure, then
\begin{align*}
\int_{\mathbb R^n}\big|f(t)\big|^2 dt\le K\,e^{K\left|E_1\right|\left|E_2\right|}\left\{\int_{\mathbb R^n\setminus E_1}\big|f(t)\big|^2 dt+\int_{\mathbb R^n\setminus E_2}\left|\hat {f}(\xi)\right|^2 d\xi\right\},\tag{4.1}
\end{align*}
where $K$ is a positive constant, and $\left|E_1\right|$ and $\left|E_2\right|$ denote the measures of  $E_1$ and $E_2$ , respectively.

\parindent=8mm \vspace{.1in}
In this Section, our primary interest is to establish the Nazarov's uncertainty principle for the continuous shearlet transforms in arbitrary space dimensions by employing the inequality (4.1). In this direction, we have the following main theorem.

\parindent=0mm \vspace{.1in}
{\bf Theorem 4.1.} {\it Let $\big[{\mathcal {SH}}_{\psi}f\big](a,s,t)$  be the shearlet transform of any  arbitrary function  $f\in L^2(\mathbb R^n)$, then the following uncertainty inequality holds:}
\begin{align*}
\int_{\mathbb R^n\setminus E_1}\int_{\mathbb R^{n-1}}\int_{{\mathbb R \setminus \left\{0\right\}}}\Big|{\mathcal {SH}}_{\psi}f(a,s,t)\Big|^2 \dfrac{da\,ds\,dt}{a^{n+1}}+C_{\psi}\int_{\mathbb R^n\setminus E_2}\left|\hat{f}(\xi)\right|^2d\xi\ge \dfrac{C_{\psi}\big\|f\big\|^{2}_{2}}{e^{K\left|E_1\right|\left|E_2\right|}},\tag{4.2}
\end{align*}
{\it where $C_{\psi}$ is given by (2.6), $E_1$, $E_2$ are two subsets of $\mathbb R^n$ with finite measures and $K$ is a positive constant.}

\parindent=0mm \vspace{.1in}
{\it Proof.} Since ${\mathcal {SH}}_{\psi}f(a,s,t)\in L^2(\mathbb R^n)$, whenever $f\in L^2(\mathbb R^n)$,  so we can replace the function $f$ appearing in (4.1) with  ${\mathcal {SH}}_{\psi}f(a,s,t)$ to get
\begin{align*}
&\int_{\mathbb R^n}\Big|{\mathcal {SH}}_{\psi}f(a,s,t)\Big|^2 dt\\
&\le K\,e^{K\left|E_1\right|\left|E_2\right|}\left\{\int_{\mathbb R^n\setminus E_1}\Big|{\mathcal {SH}}_{\psi}f(a,s,t)\Big|^2 dt+\int_{\mathbb R^n\setminus E_2}\left|\mathscr F\Big[{\mathcal {SH}}_{\psi}f(a,s,t)\Big](\xi)\right|^2 d\xi\right\}.\tag{4.3}
\end{align*}
By integrating (4.3), we obtain
\begin{align*}
&\int_{\mathbb R^{n}}\int_{\mathbb R^{n-1}}\int_{{\mathbb R \setminus \left\{0\right\}}}\Big|{\mathcal {SH}}_{\psi}f(a,s,t)\Big|^2 \dfrac{da\,ds\,dt}{a^{n+1}}\le K\,e^{K\left|E_1\right|\left|E_2\right|}\\
&\times\left\{\int_{\mathbb R^n\setminus E_1}\int_{\mathbb R^{n-1}}\int_{{\mathbb R \setminus \left\{0\right\}}}\Big|{\mathcal {SH}}_{\psi}f(a,s,t)\Big|^2 \dfrac{da\,ds\,dt}{a^{n+1}}+\int_{\mathbb R^2\setminus E_2}\int_{\mathbb R}\int_{{\mathbb R \setminus \left\{0\right\}}}\Big|\mathscr F\Big[{\mathcal {SH}}_{\psi}f(a,s,t)\Big](\xi)\Big|^2 \dfrac{da\,ds\,d\xi}{a^{n+1}}\right\}.
\end{align*}
Using Lemma 2.1 together with the energy preserving relation (2.9), the above inequality becomes
\begin{align*}
\int_{\mathbb R^n\setminus E_1}\int_{\mathbb R^{n-1}}\int_{{\mathbb R \setminus \left\{0\right\}}}\Big|{\mathcal {SH}}_{\psi}f(a,s,t)\Big|^2 \dfrac{da\,ds\,dt}{a^{n+1}}&+\int_{\mathbb R^n\setminus E_2}\int_{\mathbb R^{n-1}}\int_{{\mathbb R \setminus \left\{0\right\}}}\left|\hat{f}(\xi)~\overline{\hat{\psi}( M_{sa}\xi)}\right|^2\,\dfrac{d\xi\,ds\,da}{a^{\frac{n^2-n+1}{n}}}\\
&\qquad\qquad\qquad\qquad\ge\dfrac{C_{\psi}\big\|f\big\|^{2}_{2}}{Ke^{K\left|E_1\right|\left|E_2\right|}},
\end{align*}
which further implies
\begin{align*}
\int_{\mathbb R^n\setminus E_1}\int_{\mathbb R^{n-1}}\int_{{\mathbb R \setminus \left\{0\right\}}}\Big|{\mathcal {SH}}_{\psi}f(a,s,t)\Big|^2 \dfrac{da\,ds\,dt}{a^3}&+\int_{\mathbb R^n\setminus E_2}\left|\hat{f}(\xi)\right|^2\left\{\int_{\mathbb R^{n-1}}\int_{{\mathbb R \setminus \left\{0\right\}}}\dfrac{\Big|\hat{\psi}( M_{sa}\xi)\Big|^2}{a^{\frac{n^2-n+1}{n}}}da\,ds \right\}\,d\xi\\
&\qquad\qquad\qquad\qquad\quad\ge \dfrac{C_{\psi}\big\|f\big\|^{2}_{2}}{Ke^{K\left|E_1\right|\left|E_2\right|}}.\tag{4.4}
\end{align*}
Since $\psi\in L^2(\mathbb R^n)$ is an admissible shearlet, therefore (4.4) takes the  form
\begin{align*}
\int_{\mathbb R^n\setminus E_1}\int_{\mathbb R^{n-1}}\int_{{\mathbb R \setminus \left\{0\right\}}}\Big|{\mathcal {SH}}_{\psi}f(a,s,t)\Big|^2 \dfrac{da\,ds\,dt}{a^3}+C_{\psi}\int_{\mathbb R^n\setminus E_2}\left|\hat{f}(\xi)\right|^2d\xi\ge \dfrac{C_{\psi}\big\|f\big\|^{2}_{2}}{Ke^{K\left|E_1\right|\left|E_2\right|}},
\end{align*}
which is the desired Nazarov's uncertainty principle for the continuous shearlet transforms in arbitrary space dimensions. \qquad \fbox

\parindent=0mm \vspace{.1in}

{\it Deduction:}  As a consequence of (4.1), we can write
\begin{align*}
\int_{\mathbb R^n\setminus E_2}\left|\hat {f}(\xi)\right|^2 d\xi&\ge\dfrac{1}{K\,e^{K\left|E_1\right|\left|E_2\right|}}\int_{\mathbb R^n}\big|f(t)\big|^2 dt-\int_{\mathbb R^n\setminus E_1}\big|f(t)\big|^2 dt\qquad\quad\\
&=\dfrac{\big\|f\big\|_2^2}{Ke^{K\left|E_1\right|\left|E_2\right|}}-\int_{\mathbb R^n\setminus E_1}\big|f(t)\big|^2 dt.\tag{4.5}
\end{align*}
Using (4.5) in (4.2), the Nazrov's inequality for the shearlet transforms reduces to
\begin{align*}
\int_{\mathbb R^n\setminus E_1}\int_{\mathbb R^{n-1}}\int_{{\mathbb R \setminus \left\{0\right\}}}\Big|{\mathcal {SH}}_{\psi}f(a,s,t)\Big|^2 \dfrac{da\,ds\,dt}{a^{n+1}}+\dfrac{C_{\psi}\big\|f\big\|_2^2}{Ne^{N\left|E_1\right|\left|E_2\right|}}-C_{\psi}\int_{\mathbb R^n\setminus E_1}\big|f(t)\big|^2 dt\ge \dfrac{C_{\psi}\big\|f\big\|^{2}_{2}}{Ne^{N\left|E_1\right|\left|E_2\right|}}.
\end{align*}
Consequently, we have
\begin{align*}
\int_{\mathbb R^n\setminus E_1}\int_{\mathbb R^{n-1}}\int_{{\mathbb R \setminus \left\{0\right\}}}\Big|{\mathcal {SH}}_{\psi}f(a,s,t)\Big|^2 \dfrac{da\,ds\,dt}{a^{n+1}}-C_{\psi}\int_{\mathbb R^n\setminus E_1}\big|f(t)\big|^2 dt\ge 0.\tag{4.6}
\end{align*}
From inequality (4.6), we observe that, except for the factor $C_{\psi}$, the net concentration of the shearlet transform ${\mathcal {SH}}_{\psi}f(a,s,t)$ in $L^2\big({\mathbb R^n\setminus E_1}\times\mathbb R^{n-1}\times\mathbb R^+\big)$ is always greater than or equal to the net concentration of the signal $f$ in its natural domain $L^2\big({\mathbb R^n\setminus E_1}\big)$. Moreover, if $|E_1|=0$, then the energy preserving relation (2.9) guarantees the equality in (4.6).

\section{Local-type Uncertainty Principles for the Shearlet Transforms}

Since the classical uncertainty principle does not preclude any signal $f$ from being concentrated in a small neighbourhood of two or more widely separated points. Keeping this fact in mind, we shall derive some  local uncertainty principles for the continuous shearlet transform in arbitrary space dimensions which demonstrates that the aforementioned phenomenon can't also occur.

\parindent=0mm \vspace{.1in}

{\bf Theorem 5.1.} {\it Let $\psi$ be an admissible shearlet in $ L^2(\mathbb R^n)$. Then, for  any  $f\in L^2(\mathbb R^n)$, we have the following uncertainty inequality}
\begin{align*}
&\int_{\mathbb S}\big| t\big|^{2\alpha}\Big|{\mathcal {SH}}_{\psi}f(a,s,t)\Big|^2 d\eta\ge \dfrac{C_{\psi}}{K_{\alpha}\left|E\right|^{\alpha}}\int_{E}\left|\hat{f}(\xi)\right|^2d\xi,\quad0<\alpha<1.\tag{5.1}
\end{align*}
{\it where  $E$ is a measurable set with finite measure and $K_{\alpha}$ is a constant.}

\parindent=0mm \vspace{.1in}
{\it Proof.} For  $E\subset \mathbb R^n$ with finite measure and  $f\in L^2(\mathbb R^n)$, there exist a constant $K_{\alpha}, 0<\alpha<1$, such that  \cite{Fol}
\begin{align*}
\int_{E}\left|\hat{f}(\xi)\right|^2d\xi\le K_{\alpha}\left|E\right|^{\alpha}\Big\|\left| t\right|^{\alpha}f(t)\Big\|_2^2.\tag{5.2}
\end{align*}
Using (5.2) for the continuous shearlet transforms ${\mathcal {SH}}_{\psi}f(a,s,t)$, we obtain
\begin{align*}
\int_{E}\Big|\mathscr F\Big[{\mathcal {SH}}_{\psi}f(a,s,t)\Big](\xi)\Big|^2d\xi\le K_{\alpha}\left|E\right|^{\alpha}\Big\|\left| t\right|^{\alpha}{\mathcal {SH}}_{\psi}f(a,s,t)\Big\|_2^2.\tag{5.3}
\end{align*}
For explicit expression of (5.3), we shall integrate this inequality with respect to the measure $dads/a^{n+1}$ to get
\begin{align*}
\int_{E}\int_{\mathbb R^{n-1}}\int_{{\mathbb R \setminus \left\{0\right\}}}\Big|\mathscr F\Big[{\mathcal {SH}}_{\psi}f(a,s,t)\Big](\xi)\Big|^2\dfrac{da\,ds\,d\xi}{a^{n+1}}\le K_{\alpha}\left|E\right|^{\alpha}\int_{\mathbb S}\big| t\big|^{2\alpha}\Big|{\mathcal {SH}}_{\psi}f(a,s,t)\Big|^2d\eta,
\end{align*}
which together with Lemma 2.1 gives
\begin{align*}
\int_{E}\int_{\mathbb R^{n-1}}\int_{{\mathbb R \setminus \left\{0\right\}}}\left|\hat{f}(\xi)~\overline{\hat{\psi}( M_{sa}\xi)}\right|^2\,\dfrac{d\xi\,ds\,da}{a^{\frac{n^2-n+1}{n}}}\le K_{\alpha}\left|E\right|^{\alpha}\int_{\mathbb S}\big| t\big|^{2\alpha}\Big|{\mathcal {SH}}_{\psi}f(a,s,t)\Big|^2d\eta.\tag{5.4}
\end{align*}
Since $\psi$ is an admissible shearlet, inequality (5.4) reduces to
\begin{align*}
C_{\psi}\int_{E}\left|\hat{f}(\xi)\right|^2d\xi\le K_{\alpha}\left|E\right|^{\alpha}\int_{\mathbb S}\big| t\big|^{2\alpha}\Big|{\mathcal {SH}}_{\psi}f(a,s,t)\Big|^2d\eta.
\end{align*}
Or equivalently,
\begin{align*}
 \int_{\mathbb S}\big| t\big|^{2\alpha}\Big|{\mathcal {SH}}_{\psi}f(a,s,t)\Big|^2 d\eta\ge \dfrac{C_{\psi}}{K_{\alpha}\left|E\right|^{\alpha}}\int_{E}\left|\hat{f}(\xi)\right|^2d\xi,\quad0<\alpha<1.\tag{5.5}
\end{align*}

This completes the proof of Theorem 5.1.\quad\fbox

\parindent=8mm \vspace{.1in}

Based on the Sobolev-type uncertainty inequality (3.21), we shall derive another local uncertainty principle for the continuous shearlet transform in arbitrary space dimensions.

\parindent=0mm \vspace{.1in}

{\bf Theorem 5.2.} {\it Let   $\psi\in L^2(\mathbb R^n)$ be an admissible shearlet with $C_{\psi}=1$. Then, for  arbitrary function $f\in \mathbb H(\mathbb R^n)\cap \mathbb W_{1}^1(\mathbb R^n)$, we have}
\begin{align*}
\left(\int_{\mathbb S} \left|t\right|^2\Big|{\mathcal {SH}}_{\psi}f(a,s,t)\Big|^2d\eta\right)\ge \Bigg\{\dfrac{2}{\big\|\nabla f\big\|_{2}}\exp\left(\dfrac{{\Gamma^{\prime}}(n/2)}{\Gamma(n/2)}\right)\big\|f\big\|^{3}_{2}-\big\|f\big\|^{2}_{2}\Bigg\}.\tag{5.6}
\end{align*}
{\it provided the L.H.S of (5.6) is defined.}

\parindent=0mm \vspace{.1in}
{\it Proof.} For  $C_{\psi}=1$, we infer from (3.21) that
\begin{align*}
\left(\dfrac{{\Gamma^{\prime}}(n/2)}{\Gamma(n/2)}\right)\big\|f\big\|^{2}_{2}\le\int_{\mathbb S} \Big|{\mathcal {SH}}_{\psi}f(a,s,t)\Big|^2\ln\left(\dfrac{1+\left|t\right|^2}{2}\right)\,d\eta+\int_{\mathbb R^{n}}{\ln{\left|\xi\right|} \left|\hat{f}(\xi)\right|^2}d{\xi}.\tag{5.7}
\end{align*}
Using Jensen's inequality in (5.7), we can deduce that
\begin{align*}
\left(\dfrac{{\Gamma^{\prime}}(n/2)}{\Gamma(n/2)}\right)\le\ln\int_{\mathbb S} \dfrac{\big|{\mathcal {SH}}_{\psi}f(a,s,t)\big|^2}{\big\|f\big\|^{2}_{2}}\left(\dfrac{1+\left|t\right|^2}{2}\right)\,d\eta+\dfrac{1}{2}\int_{\mathbb R^{n}}{\ln{\left|\xi\right|^2} \dfrac{\left|\hat{f}(\xi)\right|^2}{\big\|f\big\|^{2}_{2}}}d{\xi}.\tag{5.8}
\end{align*}
To obtain a fruitful estimate of the second integral of (5.8), we set
\begin{align*}
d\rho=\dfrac{\left|\hat{f}(\xi)\right|^2}{\big\|f\big\|^{2}_{2}}d{\xi},\quad {\text {so that}}\quad\int_{\mathbb R^{n}}d\rho=1.\tag{5.9}
\end{align*}
Again by employing the Jensen's inequality, we obtain
\begin{align*}
\int_{\mathbb R^{n}}{\ln{\left|\xi\right|^2} \left|\hat{f}(\xi)\right|^2}d{\xi}&=\big\|f\big\|^{2}_{2}\int_{\mathbb R^{n}}\ln{\left|\xi\right|^2}d\rho\\
&\le\big\|f\big\|^{2}_{2}\,\ln\left\{\int_{\mathbb R^{n}}\left|\xi\right|^2 d\rho\right\}\\
&=\big\|f\big\|^{2}_{2}\ln \left\{ \int_{\mathbb R^{n}} \left|\xi\right|^2\dfrac{\left|\hat{f}(\xi)\right|^2}{\big\|f\big\|^{2}_{2}}d{\xi}\right\}\\
&=\big\|f\big\|^{2}_{2}\ln \left\{\dfrac{1}{\big\|f\big\|^{2}_{2}} \int_{\mathbb R^{n}}\big|\nabla f(t)\big|^{2}dt\right\}.\tag{5.10}
\end{align*}
Using the expression (5.10) in (5.9), we have
\begin{align*}
\left(\dfrac{{\Gamma^{\prime}}(n/2)}{\Gamma(n/2)}\right)\le\ln\int_{\mathbb S} \dfrac{\big|{\mathcal {SH}}_{\psi}f(a,s,t)\big|^2}{\big\|f\big\|^{2}_{2}}\left(\dfrac{1+\left|t\right|^2}{2}\right)\,d\eta+\dfrac{1}{2}\ln \left\{\dfrac{1}{\big\|f\big\|^{2}_{2}} \int_{\mathbb R^{n}}\big|\nabla f(t)\big|^{2}dt\right\}\\
=\ln \left\{\dfrac{1}{2\big\|f\big\|^{3}_{2}} \left\{\int_{\mathbb S} \Big|{\mathcal {SH}}_{\psi}f(a,s,t)\Big|^2\left(1+\left|t\right|^2\right)\,d\eta\right\} \left\{\int_{\mathbb R^{n}}\big|\nabla f(t)\big|^{2}dt\right\}^{1/2}\right\}.\tag{5.11}
\end{align*}

Expression (5.11) can be rewritten in a lucid manner as

\begin{align*}
\left\{\int_{\mathbb S} \Big|{\mathcal {SH}}_{\psi}f(a,s,t)\Big|^2\left(1+\left|t\right|^2\right)\,d\eta\right\} \left\{\int_{\mathbb R^{n}}\big|\nabla f(t)\big|^{2}dt\right\}^{1/2}\ge 2\exp\left\{\dfrac{{\Gamma^{\prime}}(n/2)}{\Gamma(n/2)}\right\}\big\|f\big\|^{3}_{2}.\tag{5.12}
\end{align*}

Applying the energy preserving relation (2.9) with $C_{\psi}=1$, we get

\begin{align*}
\left\{\int_{\mathbb S} \left|t\right|^2\Big|{\mathcal {SH}}_{\psi}f(a,s,t)\Big|^2d\eta\right\} \left\{\int_{\mathbb R^{n}}\big|\nabla f(t)\big|^{2}dt\right\}^{1/2}\ge 2\exp\left(\dfrac{{\Gamma^{\prime}}(n/2)}{\Gamma(n/2)}\right)\big\|f\big\|^{3}_{2}-\big\|f\big\|^{2}_{2}\big\|\nabla f\big\|_{2},
\end{align*}
which upon simplification gives the desired inequality
\begin{align*}
\left\{\int_{\mathbb S} \left|t\right|^2\Big|{\mathcal {SH}}_{\psi}f(a,s,t)\Big|^2d\eta\right\}\ge \Bigg\{\dfrac{2}{\big\|\nabla f\big\|_{2}}\exp\left(\dfrac{{\Gamma^{\prime}}(n/2)}{\Gamma(n/2)}\right)\big\|f\big\|^{3}_{2}-\big\|f\big\|^{2}_{2}\Bigg\}.
\end{align*}

This completes the proof of the theorem.\quad\fbox

\end{document}